\documentclass{IEEEtran4PSCC}

\ifCLASSINFOpdf
   \usepackage[pdftex]{graphicx}
\else
   \usepackage[dvips]{graphicx}
\fi

\usepackage[cmex10]{amsmath}
\usepackage{amsfonts,amssymb}
\usepackage{algorithmic}
\usepackage{threeparttable}
\usepackage{makecell}
\usepackage{multirow}
\usepackage{ifthen}

\newboolean{no_sticky_notes}
\setboolean{no_sticky_notes}{true}

\newcommand {\paren}[1]{\left ( #1 \right )}
\newcommand {\bparen}[1]{\left [ #1 \right ]}
\newcommand {\cparen}[1]{\left \{ #1 \right \}}
\newcommand {\vparen}[1]{\left \vert #1 \right \vert}

\newcommand{\calG}{\mathcal{G}}
\newcommand{\frakC}{\mathfrak{C}}

\newcommand{\intG}{\mathcal{G}^\mathbb{Z}}
\newcommand{\contG}{\mathcal{G}^\mathbb{R}}
\newcommand{\calS}{\mathcal{S}}
\newcommand{\calL}{\mathcal{L}}

\newcommand{\exstL}{\mathcal{L}^\dagger}
\newcommand{\newL}{\mathcal{L}^\star}

\newcommand{\llL}{\mathcal{L}^{\circ}}
\newcommand{\calD}{\mathcal{D}}
\newcommand{\calN}{\mathcal{N}}
\newcommand{\calK}{\mathcal{K}}
\newcommand{\calT}{\mathcal{T}}
\newcommand{\calB}{\mathcal{B}}
\newcommand{\calX}{\mathcal{X}}

\newcommand{\operCostOm}{C^{\text{op}}_\omega}

\newcommand{\EOm}[1]{\mathbb{E}_\zeta\bparen{#1}}
\newcommand{\invCost}{C^{\text{inv}}}

\newcommand{\xG}{x^{\calG}}
\newcommand{\xD}{x^{\calD}}
\newcommand{\xL}{x^{\calL}}
\newcommand{\xS}{x^{\calS}}

\newcommand{\XG}{X^{\calG}}
\newcommand{\XS}{X^{\calS}}

\newcommand{\xbar}{\overline{x}}

\newcommand{\pG}{p^\calG}
\newcommand{\pS}{p^\calS}
\newcommand{\pD}{p^\calD}
\newcommand{\pDtier}{p^{\calD\calK}}
\newcommand{\pSCh}{p^{\calS\text{-ch}}}
\newcommand{\pSDch}{p^{\calS\text{-dch}}}
\newcommand{\pShed}{p^{\text{sh}}}

\newcommand{\fL}{f}

\newcommand{\bg}{{b,g}}
\newcommand{\bgtom}{{b,g,t,\omega}}

\newcommand{\bstom}{{b,s,t,\omega}}
\newcommand{\bdtom}{{b,d,t,\omega}}
\newcommand{\bdktom}{{b,d,k,t,\omega}}

\newcommand{\ltom}{{\ell,t,\omega}}

\begin{document}

\title{Nodal Capacity Expansion Planning with Flexible Large-Scale Load Siting}

\author{
\IEEEauthorblockN{Tomas Valencia Zuluaga\\ Jean-Paul Watson}
\IEEEauthorblockA{Center for Applied Scientific Computing \\
  Cyber \& Infrastructure Resilience \\
Lawrence Livermore National Laboratory\\
Livermore, CA, USA\\
\{tvalenciaz, watson61\}@llnl.gov}
\and
\IEEEauthorblockN{Simon Pang}
\IEEEauthorblockA{Materials Science Division\\
Lawrence Livermore National Laboratory\\
Livermore, CA, USA\\
pang6@llnl.gov
}
}

\maketitle

\begin{abstract}
  We propose explicitly incorporating large-scale load siting into a stochastic nodal power system capacity expansion planning model that concurrently co-optimizes generation, transmission, and storage expansion. The potential operational flexibility of some of these large loads is also taken into account by considering them as consisting of a set of tranches with different reliability requirements, which are modeled as a constraint on expected served energy across operational scenarios.
We implement our model as a two-stage stochastic mixed-integer optimization problem with cross-scenario expectation constraints.
To overcome the challenge of scalability, we build upon existing work to implement this model on a high performance computing platform and exploit scenario parallelization using an augmented Progressive Hedging Algorithm. The algorithm is implemented using the bounding features of mpisppy, which have shown to provide satisfactory provable optimality gaps despite the absence of theoretical guarantees of convergence. 
We test our approach and assess the value of this proactive planning framework on total system cost and reliability metrics using realistic testcases geographically assigned to San Diego and South Carolina, with datacenter and direct air capture facilities as large loads.
\end{abstract}

\begin{IEEEkeywords}
Capacity Expansion Planning, Decomposition, Mixed-Integer Programming, Stochastic Optimization, Large Load Siting.
\end{IEEEkeywords}

\thanksto{\noindent 
This work was performed under the auspices of the U.S. Department of Energy by Lawrence Livermore National Laboratory under Contract DE-AC52-07NA27344 and was funded by LDRD project 25-SI-007.
We thank our colleagues Minda Monteagudo and Matthew Signorotti for providing the load and generation timeseries for our testcases, and Alvina Aui, Wenqin Li, and Nathan Ellebracht for providing the generation emissions data, as well as all technoeconomic data for direct air capture facilities used in all tests.
We thank Gurobi for providing the academic license used to obtain the results presented in this work.
}

\section{Introduction}
For over a decade, the electricity grid has been undergoing a transformation driven by the penetration of renewable, intermittent, and decentralized generation. While we continue to address the security and reliability difficulties posed by this transformation, the proliferation of large-scale electricity consumers, driven for example by datacenters as well as proposed direct air carbon capture (DAC) facilities, is further complicating an already taxing environment for grid planning and is a topic of increasing interest across academia, industry, and policy makers \cite{johnson_implications_2025,pge_data_2025,arwa_impact_2025,north_american_electric_reliability_corporation_large_2024}.

In conventional planning methodologies, and consequently in existing Capacity Expansion Planning (CEP) software tools, it is customary to consider electricity load as a planning parameter, over which the planner has no control. The scale of the foreseen large-load increases challenges this view and pushes planners to take a more proactive role, seeking to influence load siting decisions as well as the possible co-location of on-site generation in addition to grid investments in generation, storage, and transmission.
This approach takes advantage of the fact that subject to some constraints like access to low-latency areas for datacenters and geological suitability for carbon capture facilities, siting decisions can have a significant degree of flexibility.

Besides \textit{siting} flexibility, large loads can also provide value to the grid through \textit{operational} flexibility. Demand Response (DR) programs exemplify this fact.
Several CEP models that consider DR as a resource have been proposed before (for instance, \cite{maranon-ledesma_analyzing_2019} considers generic DR programs, \cite{ramirez_co-optimization_2016} focuses on electric vehicles, \cite{barani_residential_2024} on  residential behind-the-meter load control), but these models consider aggregations of numerous small loads distributed throughout the systems, not large concentrated loads.
In the context of emerging large loads, \cite{arwa_impact_2025} proposes co-optimization of generation, storage and DAC siting and sizing (but no datacenters), in a similar fashion to what we propose here, with a more detailed modeling of DAC operation.
The authors in \cite{shao_stochastic_2025} consider co-optimization of CEP and datacenter siting and sizing with detailed modeling of datacenter operations, but do not consider DAC.
The two references above consider aggregate zonal systems with four or less representative days; in contrast, we consider a nodal representation of the system and propose a methodology aimed at larger numbers of representative scenarios.

We do not model the operational details of DAC or datacenters but instead propose a reliability-constrained CEP model, which is also a very mature field of research. Our approach is close to that of \cite{zampara_capacity_2025}, where a limit on expected unserved energy is formulated as an expectation constraint across operational scenarios, and then dualized to enable scenario decomposition. In \cite{zampara_capacity_2025}, this is only done for a single system-wide constraint instead of separately for each tranche of large load in the system like we do in this work.
In \cite{jacobson_computationally_2024}, the authors introduce budgeting variables to the first stage to separate scenario-coupling constraints. This same approach is extended in \cite{pecci_regularized_2025} for a larger-scale approach decoupling scenarios and time. In both of these cases, the constraints targeted include policy constraints like unserved energy and carbon emissions.
Reliability constraints have been included in CEP models in numerous other ways, including chance constraints \cite{rashidaee_linear_2018}, risk-averse formulations \cite{da_costa_reliability-constrained_2021}, and robust and multi-level approaches~\cite{dehghan_reliability-constrained_2016}. We refer the interested reader to \cite{da_costa_reliability-constrained_2021,cho_recent_2022} for an overview of existing approaches in the literature.
The idea of breaking load into tranches with differentiated reliability is well-known, and is implemented routinely at the consumer level through various mechanisms, including physically separated circuits (like in hospitals) and smart devices like controllable thermostats.
Although we do not go into the details of ground-level implementation, the framework we propose here is inspired by the same concept.

It has been found that improving the spatial resolution of CEP models improves the quality of the investment plans obtained~\cite{jacobson_quantifying_2024}.
Higher spatial resolution and  reliability constraints can exacerbate existing computational challenges. Numerous authors have worked on leveraging decomposition and High Performance Computing to overcome these challenges. Benders and cut-based approaches are the most commonly used techniques to this end (including the aforementioned \cite{da_costa_reliability-constrained_2021,jacobson_computationally_2024,pecci_regularized_2025,zampara_capacity_2025}). In this work, we implement 
an augmented Progressive Hedging approach, which has also been successfully used in CEP in previous work~\cite{munoz_scalable_2015,valencia_zuluaga_parallel_2024}.

In this paper, we address the challenges mentioned above by proposing a \textit{nodal} power grid expansion planning model that explicitly considers sizing and siting large loads as part of the CEP co-optimization model, and considers their operational flexibility. The concrete contributions of our work are listed below.

\subsubsection*{Contributions}
\begin{itemize}
\item We introduce a two-stage nodal mixed-integer stochastic optimization model for capacity expansion planning that explicitly considers large loads, like DAC and datacenters, as part of the planning process, in addition to co-optimizing  generation, transmission, and storage expansion.
\item We propose a framework to represent the flexibility that large consumers can provide to the grid by breaking each large load into reliability tranches and adding expectation constraints to guarantee that the model captures both the limitations and the incentives of utilizing flexible large-load resources.
\item We develop an augmented Progressive Hedging Algorithm to decompose the expectation constraints added so that a scenario decomposition can be implemented in a parallel computing cluster, and implement it with the open-source library \textit{mpisppy}.
\item We illustrate through examples the value and promise of our model and solution technique.
\end{itemize}

\subsubsection*{Structure of this paper}
Section \ref{sec:model} provides a description of the base CEP model and the modifications proposed to explicitly incorporate large load siting.
Section \ref{sec:solution} describes the solution methodology adopted to decompose the problem by scenario and corresponding solution in a parallel computing platform.
Section \ref{sec:results} presents and discusses some numerical tests performed; we close with some remarks in Section \ref{sec:conclusion}.

\section{Model}
\label{sec:model}

\subsection{Base model}
\subsubsection{Model formulation}

The basis for our CEP model is that of \cite{musselman_climate-resilient_2025,valencia_zuluaga_parallel_2024}. It consists of a two-stage stochastic program that co-optimizes generation, transmission, and storage investments in the first stage and solves a multi-period optimal operation problem in the second stage. 
Second-stage subproblems, i.e., representative days, differ from each other in the hourly nodal demand and generation availability.
This general framework is not uncommon in the literature of CEP \cite{munoz_scalable_2015, go_assessing_2016, zampara_capacity_2025}.

Inverter-based generation and storage investments are modeled as continuous variables, turbine-based generation investments as integer variables (number of installed units), 
and transmission investments as binary  (build or no-build for candidates) variables. 
All investment costs are assumed to be linear in the size of the installed capacity. 
Storage levels are assumed to be cyclic.
See \cite{musselman_climate-resilient_2025} for a detailed discussion of the base model.

Next we provide a full, but compact description of the existing model with the minor changes made to incorporate load siting. 
In sections \ref{sec:model_load_flex} and \ref{sec:model_load_incentive} we present more significant model extensions to handle the flexibility of these resources.
A full list of the symbols used in this model is provided in Appendix \ref{sec:nomenclature}.

\subsubsection*{Objective}
The objective is to minimize the total cost of expanding, maintaining, and operating the system, which we break into an investment cost $\invCost$ and an expected operation cost $\EOm{\operCostOm}$. The random variable $\zeta$ is used to represent in compact form all  scenario-dependent data.

\begin{subequations}
\begin{multline}
\invCost = 
\sum_{\ell\in\newL}C^{\calL f}_\ell\xL_\ell + 
+\sum_{b\in\calB}\left( \right. \\
  \left.\sum_{g\in\calG}C^{\calG f}_g P^\calG_g\xG_\bg
  +\sum_{s\in\calS}C^{\calS f}_s\xS_{b,s}                         
  +\sum_{d\in\calD}C^{\calD f}_s\xD_{b,d}                         
\right)
\end{multline}
\begin{multline}
\EOm{\operCostOm} = \sum_{\omega\in\Omega}\pi_\omega \operCostOm = 
\sum_{\omega\in\Omega}\left(\pi_\omega\cdot 365\sum_{t\in\calT}\sum_{b\in\calB}\tau\cdot\Bigg[
 \right.\\
  C^{\text{sh}}\pShed_{b,t,\omega} +
  \\
\left.
\left.  
  +\sum_{g\in\calG} C^{\calG v}_g\pG_\bgtom  
  +\sum_{s\in\calS} C^{\calS v}_s \pSDch_\bstom
  +\sum_{d\in\calD}C^{\calD v}_d\pD_\bdtom
  \right]  
\right)
\end{multline}
\label{eq:cep_obj_EF}
\end{subequations}

The operation cost during each representative day for a scenario $\omega$ is multiplied by 365 so that all costs are in \$/year.


\subsubsection*{Decision variables}
\paragraph{Investment variables}
The domains of the investment variables are those of (\ref{eq:cep_domains_invest}).
For large-scale loads, we only consider the case where investment variables must be integer multiples of a unit size. The model could easily be extended to also consider continuous load types.

\begin{subequations}
    \begin{align}
        \xG_\bg \in \mathbb{Z}^+ ~ \forall b\in\calB, \forall g\in \intG \\
        \xG_\bg \in \mathbb{R}^+ ~ \forall b\in\calB, \forall g\in \contG \\
        \xS_{b,s} \in \mathbb{R}^+ ~ \forall b\in\calB, \forall s\in \calS \\
        \xD_{b,d} \in \mathbb{Z}^+ ~ \forall b\in\calB, \forall d\in \calD \\
        \xL_\ell \in \cparen{0,1} ~ \forall \ell \in \newL
    \end{align}
    \label{eq:cep_domains_invest}
\end{subequations}

\paragraph{Operation variables}
All operation variables are continuous (we do not include unit commitment in the second-stage problem). All operation variables except branch flows are non-negative.


\subsubsection*{Constraints}

\paragraph{Construction constraints}
We consider construction limits for new generation and storage resources.
Limits are considered for each generation, storage, and load type at each bus (\ref{eq:cep_construction_limits}).
Limits across buses or types were omitted here but can easily be added to the model.

\begin{subequations}
    \begin{align}
        P^\calG_g\xG_\bg+\XG_\bg \le K^{\calB\calG}_\bg ~ \forall b\in\calB,g\in\calG \label{eq:cep_pmax_bg}\\
        \xS_{b,s}+\XS_{b,s} \le K^{\calB\calS}_{b,s} ~ \forall b\in\calB,s\in\calS \label{eq:cep_smax_bs} \\
        \xD_{b,d} \le K^{\calB\calD}_{b,d} ~ \forall b\in\calB,d\in\calD \label{eq:cep_dmax_bd}
    \end{align}
    \label{eq:cep_construction_limits}
\end{subequations}

\paragraph{Physical limits of generation and storage}
The output of each generator cannot exceed the available power, which is dictated by the installed capacity and the availability of intermittent resources (\ref{eq:cep_pmax_genav}).
Similarly, the charge and discharge of storage facilities cannot exceed the installed power conversion capacity (\ref{eq:cep_smax_ch}),(\ref{eq:cep_smax_dch}), and the availability of storage must be respected (\ref{eq:cep_smax_ps}).
We do not include ramping constraints or unit commitment in this model (cf. Section \ref{sec:model_limitations}), but to avoid unrealistic ramping of baseload nuclear generation, we add a minimum power output for that generation type only \eqref{eq:pmin_nuke}.

\begin{subequations}
    \begin{align}
    ~ \forall b\in\calB,t\in\calT,\omega\in\Omega ~: \nonumber \\
        \pG_\bgtom \le \alpha_\bgtom\paren{\XG_\bg + P^\calG_g\xG_\bg} ~ \forall g\in\calG\label{eq:cep_pmax_genav} \\
        \pG_\bgtom \ge P^{\calG-\text{min}}\paren{\XG_\bg + P^\calG_g\xG_\bg} ~ \forall g\in\calG^{\calN}\label{eq:pmin_nuke} \\
        \pSCh_\bstom \le \xS_{b,s} + \XS_{b,s} ~  \forall s\in\calS\label{eq:cep_smax_ch} \\
        \pSDch_\bstom \le \xS_{b,s} + \XS_{b,s} ~  \forall s\in\calS \label{eq:cep_smax_dch} \\
        \pS_\bstom \le U^{\calS}_s\paren{\xS_{b,s} + \XS_{b,s}} ~ \forall s\in\calS \label{eq:cep_smax_ps}
    \end{align}
    \label{eq:cep_g_s_limits}
\end{subequations}

\paragraph{Energy storage}
The change in energy storage level at each facility is driven by its charge and discharge (\ref{eq:cep_stor_dyn_t}). To avoid end-of-horizon effects, we use the last period of the horizon as the initial storage state (\ref{eq:cep_stor_dyn_t0}).
A constraint to impede simultaneous charging and discharging is not included in the model, as found to be unnecessary for planning purposes in previous work~\cite{go_assessing_2016,valencia_zuluaga_parallel_2024} in the absence of policies that could incentivize spurious losses.

\begin{subequations}
    \begin{align}
    ~\forall b\in\calB,s\in\calS,\omega\in\Omega,t\in\calT\setminus\cparen{\vparen{T}} : \nonumber\\
        \pS_\bstom = \pS_{b,s,t-1,\omega} + \tau\paren{\eta^{\calS\text{-ch}}_s\pSCh_\bstom - \pSDch_\bstom} \label{eq:cep_stor_dyn_t} \\
        \pS_{b,s,0,\omega} = \pS_{b,s,\vparen{\calT}-1,\omega} + \tau\paren{\eta^{\calS\text{-ch}}_s\pSCh_\bstom - \pSDch_\bstom} \label{eq:cep_stor_dyn_t0}
    \end{align}
    \label{eq:cep_stor_dyn}
\end{subequations}

\paragraph{Load shedding}
Shedded load cannot exceed demand (\ref{eq:cep_shed_load}). Note this excludes the large loads that are part of the planning process.

\begin{align}
    \pShed_{b,t,\omega} \le D_{b,t,\omega} ~ \forall b\in\calB,t\in\calT,\omega\in\Omega \label{eq:cep_shed_load}
\end{align}

\paragraph{Power balance and power flow}
We consider the standard $b\theta$ formulation of power flow (\ref{eq:power_flow_existing}), with its corresponding big-M version for candidate transmission lines (\ref{eq:power_flow_new}). 
The energy balance is ensured at each individual bus (\ref{eq:cep_nodal_balance}).

\begin{subequations}
\begin{equation}    
\fL_\ltom = b_\ell (\theta_{o(\ell)}-\theta_{d(\ell)}) ~ \forall \ell\in\exstL, t\in\calT, \omega\in\Omega \label{eq:power_flow_existing}
\end{equation}
\begin{multline}    
-M(1-\xL_\ell)\le\fL_\ltom - b_\ell (\theta_{o(\ell)}-\theta_{d(\ell)}) \\
\le M(1-\xL_\ell) ~ \forall \ell\in\newL, t\in\calT, \omega\in\Omega \label{eq:power_flow_new}
\end{multline}    
\end{subequations}

\begin{multline}
\forall b\in\calB,t\in\calT, \omega\in\Omega: \\
    \sum_{g\in\calG} \pG_\bgtom + \sum_{s\in\calS}\paren{\eta^{\calS\text{-dch}}_s\pSDch_\bstom - \pSCh_\bstom}\\
    ~ ~ -\sum_{\ell\in\calL:o(\ell)=b}\fL_\ltom+\sum_{\ell\in\llL:d(\ell)=b}\fL_\ltom \\    
    ~ ~ +\pShed_{b,t,\omega} - \sum_{d\in\calD}\pD_\bdtom
    = D_{b,t,\omega} 
    \label{eq:cep_nodal_balance}
\end{multline}

\paragraph{Transmission limits}
Limits for existing branches are enforced by (\ref{eq:cep_max_br_flow_exst_ll}), while (\ref{eq:cep_max_br_flow_new_ll}) ensures them for built branches, as well as that unbuilt branches have no flow.

\begin{subequations}
\begin{align}
    -F_\ell\le \fL_\ltom \le F_\ell ~ \forall \ell\in \exstL,t\in\calT, \omega\in\Omega  \label{eq:cep_max_br_flow_exst_ll} \\    
    -F_\ell\xL_\ell\le \fL_\ltom \le F_\ell\xL_\ell ~ \forall \ell\in \newL, t\in\calT, \omega\in\Omega  \label{eq:cep_max_br_flow_new_ll}   
\end{align}
    \label{eq:cep_max_br_flow}
\end{subequations}

\subsubsection{Limitations of the model}
\label{sec:model_limitations}

We do not include in this model some features that are important for power system flexibility, like unit commitment (UC) or generation ramping constraints.
We omit these constraints to avoid increasing the computational complexity of our optimization problem: UC implies binary variables in the second stage, ramping constraints often require sub-hourly resolution to be binding.
A similar reasoning is followed to omit N-1 security constraints.
In the absence of such constraints, the flexibility of the system is overestimated. 
As a result, the value of considering the flexibility of large loads identified in this model is only an underestimation of their real value to the system.

\subsection{Incorporating flexibility of large loads}
\label{sec:model_load_flex}

\begin{figure}[!ht]
\centering
\includegraphics[width=2.3in]{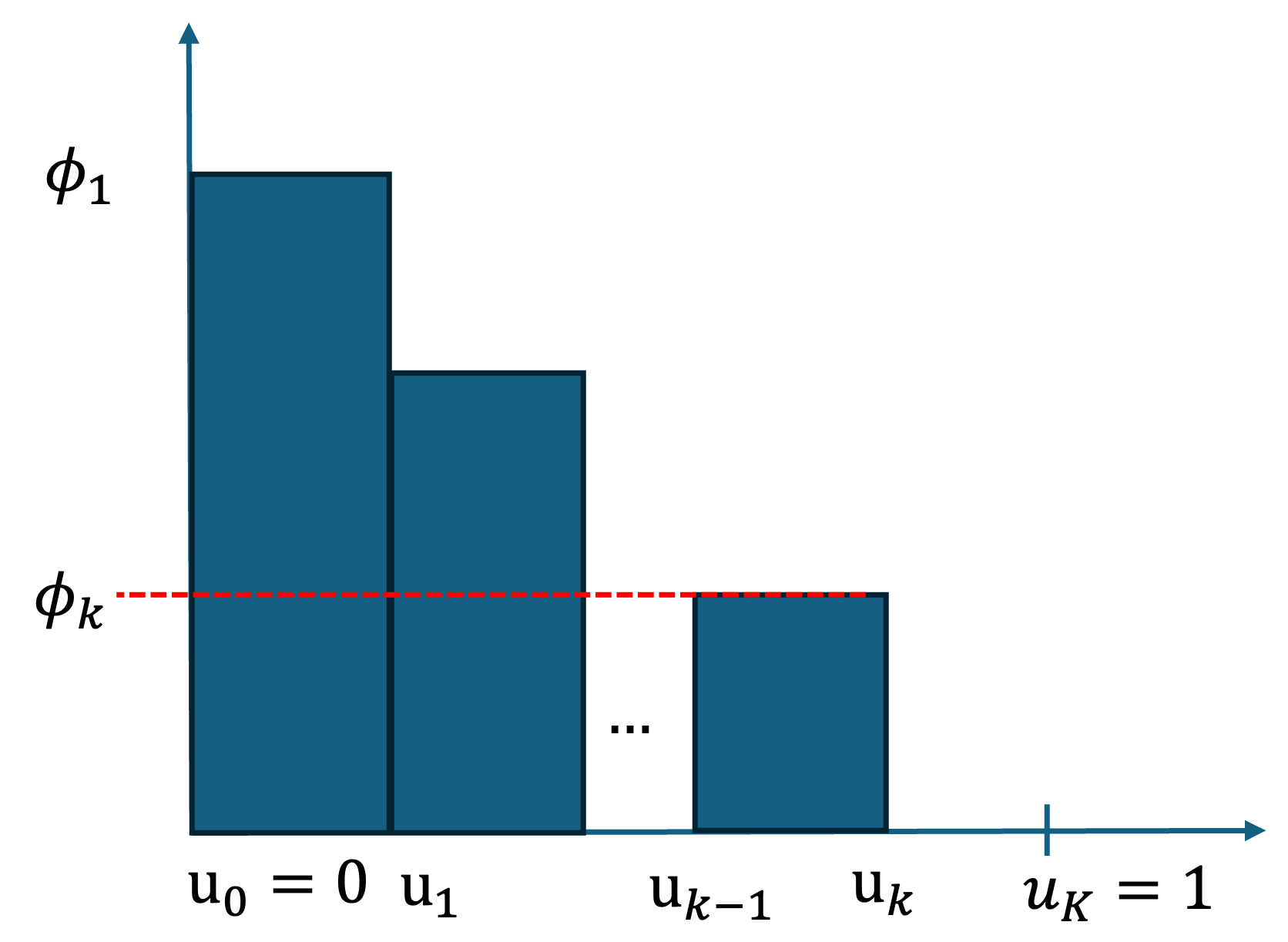}
\caption{Approach used to represent flexibility of large loads. Loads are assumed to be composed of tiers $k=1,\dots,\vparen{\calK}$, with tier $k$ consisting of a fraction $u_k-u_{k-1}$ of total load, and requiring expected service of at least $\phi_k$.}
\label{fig:tiered_model}
\end{figure}

Large loads that are included in the planning process are assumed to be composed of tiers with different reliability requirements.
Each tier is defined by the fraction of total demand that it encompasses and its required reliability level, expressed as the expected capacity factor it must have.
The scheme adopted to represent this is illustrated in Figure \ref{fig:tiered_model}. Each tier $k=1,2,\dots,\vparen{\calK}$ is defined by a pair $(u_k,\phi_k)$, with $u_k,\phi_k\in\bparen{0,1}$.
Tier $k$ consists of a fraction $u_k-u_{k-1}$ (with $u_0=0$) of the total installed capacity and must be served with a capacity factor of at least $\phi_k$.

Without loss of generality, we may sort the tiers by level of reliability and impose that $\phi_k$ be non-increasing. Note we must also have $u_{\vparen{\calK}}=1$.
Two special single-tier cases of worth are a completely inflexible load, given by $(u_1,\phi_1)=(1,1)$ and a fully flexible load, given by  $(u_1,\phi_1)=(1,0)$.
Finally, note that for simplicity we consider that all load types have the same number of tiers $\vparen{\calK}$. This is without loss of generality since we can always add dummy tiers $(0,1)$ to any load type.

\paragraph{Flexibility tier constraints}
Consistency of flexibility tiers is ensured by (\ref{eq:pDtier_max}) and (\ref{eq:pDSum}), while (\ref{eq:exp_pDtier}) ensures that the required reliability is respected.

\begin{gather}
\forall b\in\calB,d\in\calD,t\in\calT,\omega\in\Omega: \nonumber\\
    \pDtier_\bdktom \le \paren{u_{d,k}-u_{d,k-1}} \xD_{b,d} ~ ,\forall k\in\calK \label{eq:pDtier_max}\\
    \sum_{k\in\calK}\pDtier_\bdktom = \pD_\bdtom  \label{eq:pDSum} \\
    \sum_{\omega\in\Omega}\pi_\omega\sum_{t\in\calT}\tau\cdot\pDtier_\bdktom \ge \phi_{d,k} \paren{u_{d,k}-u_{d,k-1}} \vparen{\calT}\tau\ \xD_{b,d} \label{eq:exp_pDtier}
\end{gather}

\subsection{Incentivizing construction and operation of large loads}
\label{sec:model_load_incentive}

Note that because the objective is cost minimization, it would be trivially optimal to not build any large loads and thus avoid the associated costs unless the model includes some incentive for their construction and operation. We discuss next different ways in which that can be addressed here.

\subsubsection{Mandated expansion}
Perhaps the most direct and obvious way of incentivizing construction is by mandating it via a constraint like (\ref{eq:mandated_load_expansion}) below.
This could be the case when system-wide expansion plans already exist, which could come from emissions policy in the case of DAC facilities, or expansion plans of private developers in the case of datacenters.

\begin{equation}
    \sum_{b\in\calB}\xD_{b,d} \ge M_d \forall d\in \calD \label{eq:mandated_load_expansion}
\end{equation}

Note that if the type in question is fully inflexible, i.e. $(u_{d,1},\phi_{d,1})=(1,1)$, this mandate also provides a sufficient incentive for serving the load sites built.

\subsubsection{Incentive to produce output}

For flexible types, an additional incentive must be given to utilize flexible resources; otherwise, it would be trivially optimal to only utilize each tier to its minimum service level required.

\paragraph{Via expectation constraints}

We can provide such an incentive via expectation constraints imposing some bound on the output of those large loads, which could take the form of constraint (\ref{eq:expected_output_abstract}).

\begin{equation}
    \sum_{\omega\in\Omega}\pi_\omega \paren{q\cdot \pG + r\cdot\pD }\le E
    \label{eq:expected_output_abstract}
\end{equation}

The interpretation of vectors $q$ and $r$, and scalar $E$,  differs depending on the type of load.
For the case of carbon capture, $q$ could encode emission factors for each generation type, and $r$ could encode capture factors for each DAC type, so that $E$ would represent a net emissions target that may be set by policy.
For the case of datacenters, $q=0$ and $r<0$ encodes how consumed energy translates to a metric of interest, with $E<0$ and $-E$ being the minimum threshold for that metric, for example server availability, or total training core-hours of AI models.

Note that if the threshold $E$ for the metric of interest has a non-trivial value (not 0), an expectation constraint in the style of (\ref{eq:expected_output_abstract}) also provides an incentive for expansion, so that a mandate would be unnecessary.

\paragraph{Via negative costs}

It is also possible to incentivize the utilization of flexible resources by encoding the economic value of their output into the objective, by substracting some reward term $r_d$ it from its variable cost, such that  $-r_d+C^{\calD v}_d<0$.
Presumably, the resulting value of $r_d-C^{\calD v}_d$ would be in an intermediate range among the levelized cost of energy of the different generation types considered for expansion. Otherwise, the solution would trivially be to always or never serve that load type.
In our experience, even in this case, having a net negative cost with a free sizing decision $\xD_d$ leads to extreme solutions where load sites and generation capacity are built to their allowed maximum. This is a consequence of our simple linear model not taking into account the diminishing economic returns of the output $\pD_d$. 

While this could be addressed by considering a piecewise linear reward instead of a single-block $r_d$, we consider that case out of the scope of this paper. Instead, we limit models where costs are made negative to cases where we only have a load siting (not sizing) problem, i.e., when constraint (\ref{eq:mandated_load_expansion}) takes an equality form.

\section{Solution approach}
\label{sec:solution}
For instances of moderate size, the mixed-integer linear program described above can be solved in Extensive Form $(EF)$. This is, however not a scalable approach as the size of the problem increases in number of scenarios and/or power system size.

\begin{align*}
    (EF):~&\min \left( \ref{eq:cep_obj_EF} \right) \\
    &\text{s.t.} \left( \ref{eq:cep_domains_invest} \right) - \left( \ref{eq:expected_output_abstract} \right)
\end{align*}

\subsection{Relaxing expectation constraints}

Large instances will necessitate decomposition, so the expectation constraints introduced need to be reformulated. 
As is standard practice, we will dualize these constraints.
Let $\frakC$ be the union set of all expectation constraints defined in (\ref{eq:exp_pDtier}) and (\ref{eq:expected_output_abstract}), indexed by $c$. We write each constraint $c\in\frakC$ in the form (\ref{eq:abstract_exp_constr}), then 
 introduce a slack variable $\sigma_{c,\omega}$ given by (\ref{eq:abstract_exp_constr_sigma}) and associate
 a Lagrangian multiplier $\lambda_c\ge0$ with it.
We use here abstract symbols $x$ and $y_\omega$ to refer to the vectors of all first-stage and second-stage variables respectively, and $e_c$, $f_c$ and $h_c$ to refer to a scalar and vectors of coefficients of appropriate size.
 
\begin{gather}
    \sum_{\omega\in\Omega}\pi_\omega \paren{ f_c\cdot x+ h_c\cdot y_\omega} \ge e_c\label{eq:abstract_exp_constr} \\
    \sigma_{c,\omega} = e_c -  f_c\cdot x+ h_c\cdot y_\omega
    \label{eq:abstract_exp_constr_sigma}
\end{gather}

Doing this allows introducing the Lagrangian relaxation of $(EF)$, which we call $(LR)(\lambda)$.

\begin{align*}
    (LR)({\lambda}):~&\min \left( \ref{eq:cep_obj_EF} \right) + \sum_{\omega\in\Omega}\pi_\omega\sum_{c\in\frakC}\lambda_c\sigma_{c,\omega}\\
    &\text{s.t.} \left( \ref{eq:cep_domains_invest} \right) - \left( \ref{eq:pDSum} \right)
\end{align*}

Through standard duality arguments, we can conclude that $(LR)(\lambda)$ provides a lower bound of $(EF)$. Moreover, if strong duality holds, the optimal solution of $(EF)$ can be obtained 
by solving the dual problem $\max_{\lambda\ge0} LR(\lambda)$. Note that strong duality will not be guaranteed in the presence of integers or binaries, which exist in this model.

\paragraph*{Guaranteeing the feasibility of $(EF)$}
In general, problem $(EF)$ might not be feasible. Relaxing and dualizing the constraints in the manner exposed here ensures feasibility, which is convenient for computational implementations. 
A characterization of the conditions that guarantee feasibility is out of the scope of this paper, but it could be informative to consider a couple of sufficient conditions, other than the 
trivial option of adding a slack variable to each constraint \eqref{eq:expected_output_abstract}.
For reliability constraints, a sufficient condition is to add a fully dispatchable (i.e. $\alpha_{b,g,t,\omega}=1$) generator candidate at each large load location with nameplate capacity no lower than the load, as that guarantees being able to meet it with 100\% reliability (possibly at the cost of shedding load elsewhere).
For CO$_2$ emissions, a sufficient condition to guarantee feasibility is the existence of at least one site with unlimited expansion capacity for dispatchable generation and DAC facilities such that their combined emission factor is negative.

\subsection{Decomposition approach for large instances}

By construction, the second-stage variables in $(LR)$ are no longer coupled across scenarios other than by the non-anticipativity of first-stage variables $x$.
This methodology can thus be combined with existing decomposition techniques for two-stage stochastic optimization programs, like the Progressive Hedging Algorithm (PHA).

Recall that in the PHA, a copy of each first-stage variable is created for each scenario, and an equality constraint making all copies equal to each other is relaxed and dualized. 
In addition to the conventional Lagrangian term added to the objective, a so-called proximal term is added, resulting in problem $(PHA)_\omega$ for each scenario $\omega\in\Omega$.
For a more detailed description of the PHA, we refer the  reader to \cite{rockafellar_scenarios_1991,kaisermayer_progressive_2021}.

\begin{align*}
    (PHA)_\omega({\lambda}):~&\min \pi_\omega \cdot\left( C^{inv} + C^{op}_\omega 
    +\sum_{c\in\frakC}\lambda_c\sigma_{c,\omega} \right. \\
    &~ ~ \left.+ {w}\cdot x_\omega + \rho \paren{x_\omega-\xbar}^2 \right)     \\
    &\text{s.t.} \left( \ref{eq:cep_domains_invest} \right) - \left( \ref{eq:pDSum} \right) \text{ for } \omega \text{ only }
\end{align*}

We propose to solve our problem via the augmented PHA described below.

\vspace{0.5cm}
\begin{algorithmic}[1]
\STATE Initialize $w_\omega^0\gets 0~\forall\omega,\lambda_c\gets\lambda^0_c\forall c\in\frakC$
, $k \gets 0$, choose $\rho,\beta > 0$
\REPEAT
    \FOR{each scenario $\omega$}
        \STATE Solve scenario subproblem $(PHA)_\omega$to get $x_\omega^{k+1}$:
        \STATE \hspace{1em} $x_\omega^{k+1} = \arg\min_{x_\omega \in \calX_\omega} \left\{ C^{inv}+C^{op}_\omega \right.$
        \STATE \hspace{8em} $\left. w_\omega^k x_\omega + \frac{\rho}{2} \| x_\omega - \bar{x}^k \|^2 \right.$
        \STATE \hspace{8em} $\left. +\sum_{c\in\frakC}\lambda^k_c\sigma_{c,\omega} \right\}$
    \ENDFOR
    \STATE Update avg: $\bar{x}^{k+1} = \sum_\omega \pi_\omega x_\omega^{k+1}$
    \FOR{each scenario $\omega$}
        \STATE Update multipliers: $w_\omega^{k+1} = w_\omega^k + \rho (x_\omega^{k+1} - \bar{x}^{k+1})$        
    \ENDFOR
    \FOR{each expectation constraint $c\in\frakC$}
        \STATE Update avg: $\bar{\sigma}^{k+1}_c = \sum_\omega \pi_\omega \sigma_{c,\omega}^{k+1}~\forall c\in\frakC$
        \STATE Update multipliers: $\lambda_c^{k+1} = \lambda_c^k + \beta \bar{\sigma_c}$
    \ENDFOR        
    \STATE $k \leftarrow k + 1$
\UNTIL{convergence criterion}
\end{algorithmic}

\subsubsection*{About convergence guarantees}

In the convex case, standard strong duality arguments can be used to guarantee convergence of the algorithm for an adequate choice of the step sizes $\beta$ and $\rho$ (see for example \cite{zampara_capacity_2025,boyd_convex_2004}).
Because of the presence of integers and binaries in the first-stage variables, our problem is not convex and thus neither convergence nor optimality at convergence are guaranteed.
Consequently, this augmented PHA is, as the original PHA itself, only a heuristic technique on our non-convex case.
However, it can be a satisfactory heuristic if methods are leveraged to obtain valid bounds that provide
quality guarantees of the solution obtained.

\subsubsection{Implementation in \textit{mpisppy}}

We implement our model in Pyomo, and use \textit{mpisppy} to run the PHA  utilizing parallel computing for the different scenario subproblems. 
The \textit{mpisppy} package \cite{knueven_parallel_2023,woodruff_pyomompi-sppy_2025} allows implementation in a parallel computing cluster with minimal effort, as it takes care of the necessary communication within the computing cluster during the PHA update steps.
To implement the subgradient update feature that we describe here, we implement an mpisppy ``Extension'' to handle communication between the different compute nodes so that the update steps can be performed.

\paragraph{Obtaining upper and lower bounds on the solution}

A key feature of mpisppy is its ability to provide quality guarantees during the execution of the algorithm. The package allows dedicating certain groups of computing nodes, referred to as \textit{spokes}, to obtaining lower and upper bounds simultaneously as the PHA executes. Lower bounds may be obtained by solving relaxations of the problem, while upper bounds are obtained by testing feasible candidates. 
We refer the interested reader to our previous work \cite{valencia_zuluaga_parallel_2024} for more details about the implementation of a CEP model in a parallel computing cluster using mpisppy.

\paragraph{Obtaining lower bounds}
Note that for the model proposed here, the values of $\lambda$ and/or $\bar{\sigma}$ need to be communicated to the bounder spokes. Custom versions of the mpisppy \textit{Lagrangian} and \textit{reduced costs} spokes \cite{woodruff_pyomompi-sppy_2025} were developed for this purpose. 

\paragraph{Obtaining upper bounds}
In the conventional PHA, obtaining an upper bound consists in evaluating all second-stage subproblems for a candidate first-stage solution. In mpisppy, this is done by fixing the values of all non-anticipative variables to that of the evaluated candidate and solving the resulting second-stage subproblems (possibly in parallel).
To facilitate evaluating upper bounds before the values of the dual multipliers have stabilized, we write another mpisppy ``Extension'' to fix the values of all first-stage variables but continue updating the values of $\lambda$ for a number of extra iterations.
Note that in general, as discussed above, it is not guaranteed that for any given first-stage solution, all expectation constraints can be satisfied, i.e. we don't necessarily have complete recourse.
Once the values of first-stage variables is fixed, the resulting problem is an LP, and thus convexity guarantees stabilization of the values of $\lambda$. However, it is not guaranteed that this stabilization induces $\bar{\sigma}_c=0~\forall c\in\frakC$.

\section{Numerical results}
\label{sec:results}

\subsection{Test cases}

We test our model on two sandbox test systems: the standard IEEE 24-bus test case with restricted transmission capacity of \cite{go_assessing_2016}, overlayed onto the San Diego area, and 
the ACTIVSg 500-bus South Carolina test case \cite{birchfield_grid_2017}.

Both systems were extended with capacity expansion technoeconomic data and timeseries of load and generation availability downscaled from earth system models as in \cite{musselman_climate-resilient_2025}. 
The 24-bus test case is described in more detail in Appendix~\ref{sec:app_testcase}.
Some details of the 500 bus test system are reported in \cite{valencia_zuluaga_optimization_2024}.
We consider two types of large-scale loads that are included in the expansion plan: ``datacenter'' and ``DAC''.
Based on discussions with datacenter developers, two types of ``datacenter'' loads are considered: so-called ``aggregators'', representing developers that serve multiple customers, have intermediate-size sites and very little operational flexibility, and ``hyperscalers'', with very large sizes and more operational flexibility.
Table~\ref{tab:load_types} summarizes the characteristics of the types considered.
We conduct our tests on 26 representative days (12 for the 500-bus system) with different demand and generation availability.
These are sandbox testcases that are intended to represent a reasonable power system, but not necessarily a real one.
The full data will be provided upon request.

\begin{table}[!ht]
\renewcommand{\arraystretch}{1.3}
\centering
\caption{Characteristics of large-scale loads considered in tests.}
\label{tab:load_types}
\begin{tabular}{|p{0.8cm}|p{1.3cm}|p{2.2cm}|p{2.2cm}|}
\hline
Type & Unit Size & Expansion \& Utilization  Incentive & Flexibility\\
\hline\hline
DAC & 100 ktCO$_2$/y & 
\makecell[l]{
  $\bullet$ Net zero CO$_2$\\ emissions \\
  $\bullet$ 2MWh/tCO$_2$ \\captured
}
& 
\makecell[l]{
  Inflexible; Full-flex;\\
  Mid-flex: \\$u=(0.5, 0.75, 1)$,\\ $\phi=(1, 0.5, 0)$
} \\
\hline
HypScal & 800MW & 
\multirow{2}{=}{%
\makecell[l]{%
  $\bullet$ Installed capacity \\ target as in \eqref{eq:mandated_load_expansion} \\
  $\bullet$ Economic value: \\$C^{\calD v}_d=-4\text{\$/MWh}$
}} 
& 
\makecell[l]{
  Mid-flex:\\ $u=(0.5, 0.9, 1)$, \\$\phi=(1, 0.85, 0)$
} \\
\cline{1-2}\cline{4-4}
DatAgg & 100MW & 
& 
\makecell[l]{
  Mid-flex:\\ $u=(0.99, 1)$, \\$\phi=(1, 0)$
} \\
\hline
\end{tabular}
\end{table}

\subsection{Experiments}

We first conduct a series of runs solving the problem in Extensive Form (EF) on the 24-bus instances, to better illustrate the value of explicitly modeling the flexibility of large loads as proposed in this paper.
A handful of test cases are conducted to illustrate the features included in this model, in two groups as described in Table~\ref{tab:experiment_list}.
Test cases Ia, Ib and Ic test the sensitivity of the buildout to the three different levels of flexibility of DAC facilities described in Table~\ref{tab:load_types}. No datacenter additions are considered in these test cases.
Test cases IIa, IIb, and IIc consider datacenters and inflexible DAC facilities.
Case IIc corresponds to the co-optimization of large load siting and CEP that is proposed in this work.
Cases IIa and IIb assume that the new datacenters are added after the planning process is over. We assume the same buildout of case Ia and add the datacenters at the same locations obtained the optimal solution of IIc.
In case IIa, no new investments are assumed. In case IIb, we assume the datacenters are co-located with enough generation or storage to power the entire load. In particular, the hyper-scaler co-sites natural gas generation, while the datacenter aggregator adds 4h-battery storage. All cases assume a net zero CO$_2$ target.

To test the decomposition approach, we implement instances Ia, Ib, Ic and IIc on both test case systems and compare the performance of the EF and decomposition approaches.

\subsection{Results}

In all co-optimized cases, both the EF and decomposition approaches return feasible solutions, i.e., the emissions target and reliability requirements are all satisfied. 
In all cases, all renewable generation candidates are built up to the maximum allowed. This is in line with the result of previous studies that find it preferable to first decarbonize by switching to carbon-free generation technologies before resorting to direct air capture \cite{pett-ridge_roads_2023}.
Additional test cases, not reported here, where the emissions target was set at 25\% of the pre-expansion emissions (instead of zero) led to no DAC construction.
We note that a zero-emissions expectation constraint does incentivize construction of DAC facilities, even without a DAC mandate in the style of \eqref{eq:mandated_load_expansion}.

\begin{figure*}[!ht]
\centering
\begin{minipage}[c]{2.7in}
  \centering
  \includegraphics[width=2.7in]{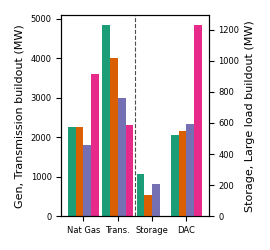}
\end{minipage}%
\begin{minipage}[c][2.5in][c]{1.5in}
  \centering  
  \includegraphics[width=1.5in,keepaspectratio]{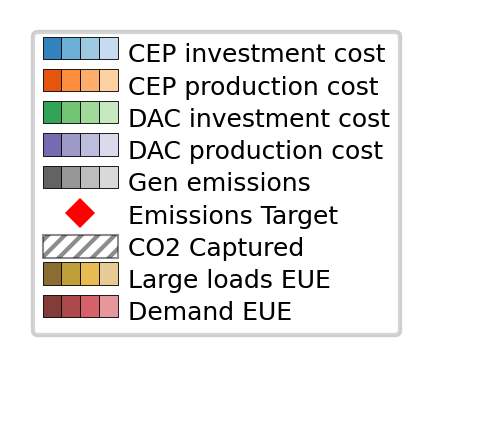}\\
  \includegraphics[width=1.5in,keepaspectratio]{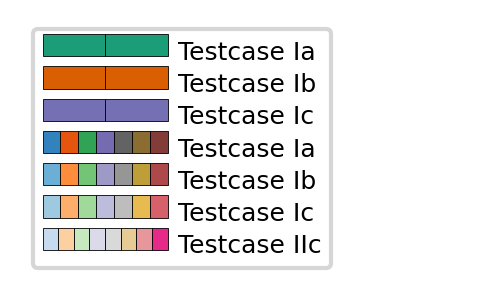}  
\end{minipage}%
\begin{minipage}[c]{2.7in}
  \centering
  \includegraphics[width=2.7in]{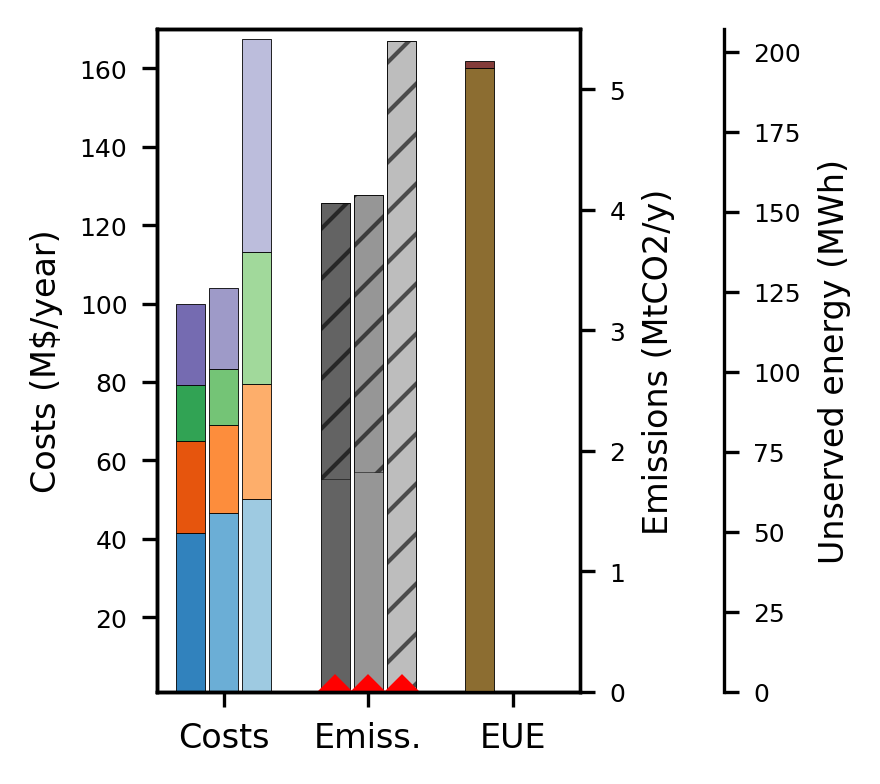}
\end{minipage}%
\caption{\textit{Left:} Resource buildout for each test case. Case IIa has by construction the same buildout as Ic. Case IIb has the same buildout plus 750MW of natural gas co-sited with the hyperscaler and 90MW of 4h battery storage for the datacenter aggregator; both are omitted from this figure.
Transmission capacity is obtained by summing the capacities of all selected candidate lines or transformers.
\textit{Right:} Total cost, CO$_2$ emissions, and expected unserved energy for test cases IIa, IIb and IIc. 
Expected Unserved Energy (EUE) refers here to shedded load ($\pShed$ in the model formulation) and expected energy shortfall for tiers served below their reliability level.
The costs exclude the penalty terms associated with slack variables of unmet expectation constraints.
}
\label{fig:case_24_results}
\end{figure*}

\begin{figure*}[!t]
\centering
\includegraphics[width=7in]{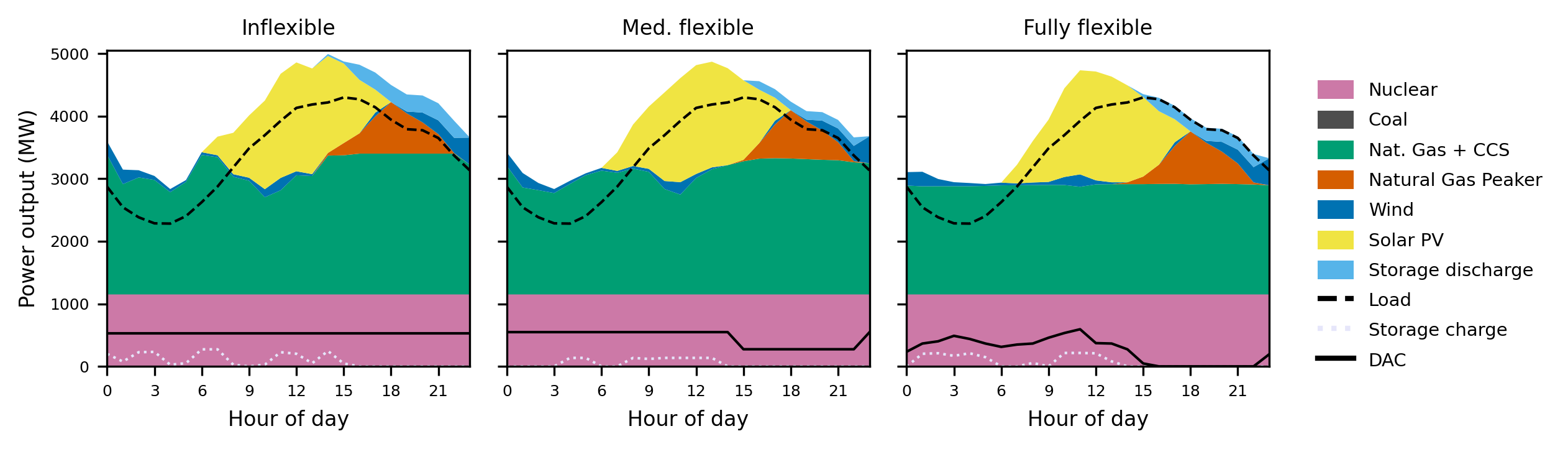}
\caption{\textit{Left:}
Operational profile of a high-load summer day with low wind availability for cases Ia (no DAC flexibility), Ib (intermediate flexibility) and Ic (full flexibility).
Additional DAC flexibility serves as peak shaving and allows operating combined cycle generation as baseload.
}
\label{fig:case_24_stackplot}
\end{figure*}

Because of space limitations, only the most illustrative results of the 24 bus test case are shown here.
Fig. \ref{fig:case_24_results}, \textit{Left} shows that as the operational flexibility of DAC facilities increases, it can replace investments in other resources used for peak shaving, like natural gas (case Ic vs Ib) or storage (Ib vs Ia). The datacenters planned in IIc are a significant addition to the test system, requiring sizable new capacity of natural gas generation and DAC.
In this example, transmission investments seem to be substituted by DAC flexibility and natural gas generation, but this result depends on congestion and may be different in other systems.
Fig.~\ref{fig:case_24_stackplot} shows a comparison of the same operational scenario (high summer load and low wind) for the three levels of DAC flexibility and illustrates the peak-shaving behavior mentioned above.
Cost differences between these cases are modest ($<2\%$) and are omitted from the plot.

Fig. \ref{fig:case_24_results}, \textit{Right} compares the costs and emissions obtained for test cases of group II.
In case IIa, it is found that failing to adequately expand the grid for the addition of such a significant load can compromise the reliability of the system, reflected in a positive expected unserved energy (EUE). Co-siting of generation alleviates the reliability issue and brings EUE to 0 (case IIb). In cases IIa and IIb, it is not possible to capture all the CO$_2$ emitted and the zero target is not met.
The co-optimization of load siting and expansion (case IIc) identifies the need for significant system additions and, as expected, allows meeting all requirements at minimal cost.

Admittedly, these are sandbox systems with a limited number of scenarios, but we believe these results provide support for continuing to develop planning tools that incorporate large loads, and pay particular attention to modeling their flexibility.

\begin{table}[!ht]
\renewcommand{\arraystretch}{1.3}
\centering
\caption{Summary of illustrative tests conducted}
\label{tab:experiment_list}
\begin{tabular}{|c|c|l|}
\hline
Test & Load type & Description \\
\hline\hline
Ia & DAC & Inflexible \\
Ib & DAC & Mid-flexible \\  
Ic & DAC & Fully flexible \\
IIa & DAC+1HypScal+1DatAgg & No extra planning \\
IIb & DAC+1HypScal+1DatAgg & Co-sited generation \\
IIc & DAC+1HypScal+1DatAgg & Proactive planning \\
\hline
\end{tabular}
\end{table}

\subsubsection*{Performance of the decomposition approach}

All tests were conducted on the \textit{dane} cluster at Lawrence Livermore National Laboratory; Gurobi 13.0 was used to solve each MILP subproblem.
We used 2 compute nodes (3 for the 500 bus test case) and a time limit of 6 hours for the entire PHA.
Table \ref{tab:computational_results} provides a summary of the computational results of the decomposition approach proposed on the different cases tested.

For the 24-bus test case, all candidate solutions obtained are near-optimal. With use of the recently developed feature of surrogate variables in \textit{mpisppy}~\cite{woodruff_pyomompi-sppy_2025}, the lower bounds were significantly tightened for cases Ia and Ib with minimal parameter tuning.
We expect the bounds of cases Ic and IIc could be tightened more with more sophisticated tuning or complementary techniques, for example implementing other spokes besides \textit{reduced costs} and Lagrangian.

For the 500-bus test cases, results are not yet satisfactory, although some recently developed \textit{mpisppy} features leveraging the use of linear model relaxations are promising.
Note however that the feasible solutions obtained with the decomposition approach are orders of magnitude better than the ones returned by the EF solver within the same time limit. It is true that the decomposition method had access to more computational resources (EF was limited to one computing node), so the comparison is not necessarily fair, but illustrates the potential of the methodology for instances where solving an EF version of the problem is not viable, which is precisely the intended use of this tool.

We expect that in order to expand to larger systems, this approach could necessitate an improved update methodology for the dual multipliers, with perhaps dynamic values for $\beta$ and $\rho$.
Another promising avenue is a hybrid approach combining the first-stage budgeting variables of \cite{jacobson_computationally_2024,pecci_regularized_2025} with our PHA methodology in the form of iterative cross-scenario cuts.
These improvements will be undertaken in future work but are outside the scope of the current paper.

\begin{table}[!ht]
\renewcommand{\arraystretch}{1.3}
\centering
\begin{threeparttable}
\caption{Performance of the decomposition approach}
\label{tab:computational_results}
\begin{tabular}{|c|cc|cc|}
\hline
 & \multicolumn{2}{c|}{\textbf{24 Buses}} & \multicolumn{2}{c|}{\textbf{500 Buses}} \\
\cline{2-5}
\textbf{Test} & Final Gap & Final Cost\tnote{a} & Final Gap & Final Cost\tnote{a,b}\,\,\\
\hline
Ia  & 2.1\% & 100\% & 89.6\% & 1.36\% \\
Ib  & 2.6\% & 99\%  & 78.4\% & 0.66\% \\
Ic  & 10.1\%    & 100\% & 90.6\%    & --\tnote{c} \\
IIc & 29.1\%   & 101\% & 96.7\% & 3.78\% \\
\hline
\end{tabular}
\begin{tablenotes}
\item[a] Normalized w.r.t.  optimal cost returned by the EF.
\item[b] For all 500-bus instances, the EF returned very poor incumbents, with $>99\%$ MIP gap.
\item[c] The EF returned no incumbent at the time limit.
\end{tablenotes}
\end{threeparttable}
\end{table}

\section{Closing remarks}
\label{sec:conclusion}

We have presented here a model that explicitly considers the siting and sizing of large loads as part of the planning process.
Our numerical results illustrate that failing to include these resources may come at the cost of meeting certain planning targets like energy not served or carbon emissions.
Moreover, explicitly including the operational flexibility of these large loads can can provide a useful signal to planners for an optimal investment plan that ensures the socially desired expansion occurs at minimum cost.

The 500-bus tests performed exhibit the computational challenges that arise with system size. The promising results obtained support continuing our work to improve the method towards an implementation of a real-sized test case with a more realistic representation of system flexibility. Future work to enable this may include more sophisticated approaches mixing the augmented PHA presented here with the recently successful cross-scenario cut-based approach with first-stage budgeting variables \cite{jacobson_computationally_2024,pecci_regularized_2025}.

\bibliographystyle{IEEEtran}
\bibliography{references_cep_dac}

\appendices

\section{List of symbols in optimization model}
\label{sec:nomenclature}
\subsection{Nomenclature}
\label{sec:cep_pha_nomenclature}

\noindent\textbf{Sets} \\
$\calB$: Set of all buses (nodes) in the network. Indexed by $b$ (or $i,m,n$ where noted).\\
$\calG$: Set of generation types, indexed by $g$, among which:\\
~~$\calG^\calN$: Short-hand for nuclear generation type.\\
~~$\intG$: Types modeled with integer variables (turbine-based).\\
~~$\contG$: Types modeled with continuous variables (inverter-based).\\
$\calS$: Set of storage types, indexed by $s$.\\
$\calL$: Set of transmission branches, indexed by $\ell$.\\
$\exstL$: Set of existing transmission branches.\\
$\newL$: Set of candidate transmission branches.\\
$\calD$: Set of large-scale load types, indexed by $d$\\
$\calK$: Set of interruptibility tiers of large loads, indexed by $k$\\
$\calT$: Set of periods in a representative day, $\cparen{0,1,\dots,\vparen{\calT}-1}$. Indexed by $t$.\\
$\Omega$: Set of scenarios in the uncertainty set of the stochastic optimization problem. Indexed by $\omega$.\\

$\intG$ and $\contG$ constitute a partition of $\calG$, i.e. $\intG\cap\contG=\emptyset$ and $\intG\cup\contG=\calG$.

\noindent\textbf{Index maps} \\
$o(\ell)$: Origin bus (also called `from' bus) of branch $\ell$\\
$d(\ell)$: Destination bus (also called `to' bus) of branch $\ell$\\

\noindent\textbf{Parameters}\\
\noindent\textit{Existing resources}\\
$\XG_{b,g}$: existing generation of type $g$ at bus $b$, in MW.\\
$\XS_{b,s}$: existing storage power conversion capacity of type $s$ at bus $b$, in MW.\\
\noindent\textit{Fixed costs}\\
$C^{(\cdot)f}_{(\cdot)}$: Annualized fixed costs (capital + fixed O\&M) for generation, storage or large load investment (in \$/MWy) \\
$C^{\calL f}_{\ell}$: Annualized fixed costs for transmission candidate (\$/y) \\
\noindent\textit{Variable costs}\\
$C^{(\cdot)v}_{(\cdot)}$: Variable costs (fuel, consumables, + variable O\&M) for generation, storage and large load investment (in \$/MWh) \\
$C^{\text{sh}}$: Cost of load shedding, in \$/MWh.\\
\noindent\textit{Operational parameters}\\
$P^\calG_g$: capacity per unit of generator type $g$, in MW.\\
$P^{\calG-min}_g$: minimum output of generator type $g$, as fraction of nameplate capacity, in MW.\\
$U^\calS_s$: Duration of storage type $s$, in hours.\\
$D_{b,t,\omega}$: demand at bus $b$ during period $t$ of scenario $\omega$, excluding large loads that are part of the planning process.\\
$\alpha_\bgtom$: fraction of generation capacity of type $g$ that is available at bus $b$ during period $t$ of scenario $\omega$.\\
$\eta^{\calS\text{-ch}}$: power conversion efficiency for  $s$ when charging.\\ 
$\eta^{\calS\text{-dch}}$: power conversion efficiency for $s$ when discharging.\\
$u_{d,k},\phi_{d,k}$: pair defining expected service level of tier $k$ of large load type $d$.\\
$b_\ell$: susceptance of branch $\ell$, in p.u.\\
$F_\ell$: transmission capacity of branch $\ell$ in MW.\\
\noindent\textit{Construction and planning parameters}\\
$K^{\calB(\cdot)}_{b,(\cdot)}$: maximum buildable capacity of specific generation, storage or large load at bus $b$, in MW.\\
\noindent\textit{Other model parameters}\\
 $\tau$: length of each period $t$ in a representative day, in h.\\
$\pi_\omega$: probability assigned to scenario $\omega$ in the stochastic optimization problem.\\
\noindent\textbf{Decision variables}\\
\noindent \textit{Investment variables}\\
$x^{(\cdot)}_{b,(\cdot)}$: new storage, generator or large load of type $(\cdot)$ at $b$.\\
$\xL_{\ell}$: binary variable indicating whether candidate $\ell$ is built.\\
\noindent \textit{Operation variables}\\
$\pG_{b,g,t,\omega}$: output of generator of type $g$ at bus $b$ during period $t$ of scenario $\omega$, in MW.\\
$\pS_\bstom$: energy stored (level) in storage facility of type $s$ at bus $b$ during period $t$ of scenario $\omega$, in MWh.\\
$\pSCh_{b,s,t,\omega}$: power input (charging) of storage facility of type $s$ at bus $b$ during period $t$ of scenario $\omega$, in MW.\\
$\pSDch_{b,s,t,\omega}$: power output (discharging) of storage facility of type $s$ at bus $b$ during period $t$ of scenario $\omega$, in MW.\\
$\pDtier_{b,d,t,k,\omega}$: power consumed by tier $k$ of large load of type $d$ at bus $b$ during period $t$ of scenario $\omega$, in MW.\\
$\pShed_{b,t,\omega}$: load shed at bus $b$ during period $t$ of scenario $\omega$, in MW (excluding large loads)\\
$\fL_\ltom$: power flow through branch $\ell$ during period $t$ of scenario $\omega$, in MW.\\
\noindent \textit{Auxiliary variables and expressions}\\
$\invCost$: Annualized total investment cost in \$/y.\\
$\operCostOm$: Operation cost of scenario $\omega$ in \$/y.\\
$\pD_{b,d,t,\omega}$: total power consumed by large load of type $d$ at bus $b$ during period $t$ of scenario $\omega$, in MW.\\
$\sigma_{(\cdot)}$: slack variable associated with an expectation constraint.\\
\noindent\textbf{Notation}\\
Unless otherwise specified, the vector is noted by omitting the corresponding index, e.g. $f = \bparen{f_\ell}_{\ell\in \calL}$.
To allow for a more compact description, $(\cdot)$ was used above as a wildcard in superscripts and subscripts and may be substituted by $\calG,g$, $\calS,s$ or $\calD,d$ as appropriate.

\section{Data for 24-bus test case}
\label{sec:app_testcase}
\paragraph{Technology types}
We consider the same basic types of \cite{go_assessing_2016}, plus natural-gas-powered combined cycle with carbon capture and sequestration. Unlike \cite{go_assessing_2016}, we don't consider power inversion capacity and storage duration distinct decision variables. Instead, we consider 2h- and 4h-duration battery storage, and 10h pumped hydro storage.
A brief overview of the technologies considered is provided in Table~\ref{tab:technology_types}.

\begin{table}[!ht]
\renewcommand{\arraystretch}{1.3}
\centering
\begin{threeparttable}
\caption{Technology types considered in the test cases}
\label{tab:technology_types}
\begin{tabular}{|c|p{4.5cm}|c|}
\hline
Label & Description & Unit capacity (MW)\\
\hline\hline
OCGT & Natural gas simple combustion turbine & 100 \\
CCGT & Natural gas combined cycle & 550 \\
CCSNG & Natural gas combined cycle with carbon capture and sequestration & 450\\
Coal & Coal-powered steam turbine & 600 \\
Nuclear & Nuclear-powered steam turbine & 800\\
Wind & Onshore wind & \tnote{*} \\
Solar PV & Solar photovoltaic & \tnote{*}  \\
PumpHyd & Pumped hydroelectric storage (10h) & \tnote{*} \\
2h BESS & 2h-duration battery energy storage  & \tnote{*} \\
4h BESS & 4h-duration battery energy storage  & \tnote{*} \\
DAC & Direct air carbon capture facility & 100\tnote{a}\\
HypScal & Hyper-scaler datacenter & 800\\
DatAgg & Aggregator datacenter & 100\\
\hline
\end{tabular}
\begin{tablenotes}
\item[*] Modeled with continuous variable; any unit capacity is only notional.
\item[a] Units of ktCO$_2$/year
\end{tablenotes}
\end{threeparttable}
\end{table}

\paragraph{Costs}
All fixed and variable cost data for generation and storage facilities is taken from \cite{musselman_climate-resilient_2025}. 

\paragraph{Expansion sites}
Solar and wind sites are only considered for expansion at the same buses as in \cite{go_assessing_2016}.
Thermal generation is considered for expansion at any bus with existing generation.
We assume battery storage can be built anywhere on the system, with a maximum of 500MW per individual bus.
Expansion sites for large loads are chosen arbitrarily aiming for diversity in location and overlap with generation candidates.
Table~\ref{tab:expansion_sites} provides an overview of the expansion sites considered.

\begin{table}[!ht]
\renewcommand{\arraystretch}{1.3}
\centering
\caption{Sites for expansion}
\label{tab:expansion_sites}
\begin{tabular}{|c|p{5.5cm}|}
\hline
Type & Buses \\
\hline\hline
Solar/Wind & 3, 5, 7, 16, 21, 23 \\
Thermal & Anywhere with existing generation \\
Battery & All buses (max 500MW at each bus) \\
Pumped Hydro & 8, 16, 22 (max 200MW each)\\
DAC & 1, 5, 8, 15, 16, 18, 23 \\
HypScal & 4, 8, 13, 19, 23 \\
DatAgg & 2, 4, 5, 14, 19 \\
\hline
\end{tabular}
\end{table}

\paragraph{Timeseries of generation and load}
Timeseries of hourly load at each bus and hourly generation availability for wind and solar sites for the different operational scenarios was obtained from downscaled reanalysis data.

All data will be provided upon request.

\end{document}